\documentclass[letterpaper,12pt,draft]{article}
\usepackage{amsmath,amssymb,amsthm, bbm}

\usepackage{graphics,graphicx}
\usepackage{pstricks,pst-node,pst-tree}
\usepackage{subfig}

\newtheorem{theorem}{Theorem}

\newtheorem{proposition}{Proposition}

\theoremstyle{definition}
\newtheorem{definition}{Definition}

\theoremstyle{remark}
\newtheorem*{remark}{Remark}
\newtheorem{note}{Note}
\newtheorem{example}{Example}

\newcommand{\ONE}{{\mathbbm{1}}}

\newcommand{\Edge}[1]{\mathcal{E}(#1)}
\newcommand{\Vertex}[1]{\mathcal{V}(#1)}

\newcommand{\Edg}{\mathcal{E}}
\newcommand{\Ver}{\mathcal{V}}
\newcommand{\Pat}{\mathcal{P}}

\makeatletter 
\def\ext@notarrow#1#2#3#4#5#6#7{%
  \mathrel{\mathop{%
    \setbox\z@\hbox{#5\displaystyle}%
    \setbox\tw@\vbox{\m@th
      \hbox{$\scriptstyle\mkern#3mu{#6}\mkern#4mu$}%
      \hbox{$\scriptstyle\mkern#3mu{#7}\mkern#4mu$}%
      \copy\z@
    }%
    \hbox to\wd\tw@{\rlap{\hbox to\wd\tw@{\hfil\raisebox{0.3ex}{$\mathbf{\scriptscriptstyle/}$}\hfil}}\unhbox\z@}}%
  \limits
    \@ifnotempty{#7}{^{\if0#1\else\mkern#1mu\fi
                       #7\if0#2\else\mkern#2mu\fi}}%
    \@ifnotempty{#6}{_{\if0#1\else\mkern#1mu\fi
                       #6\if0#2\else\mkern#2mu\fi}}}%
}
\newcommand{\xnrightarrow}[2][]{\ext@notarrow 0359\rightarrowfill@{#1}{#2}}
\newcommand{\xnleftarrow}[2][]{\ext@notarrow 3095\leftarrowfill@{#1}{#2}}
\makeatother

\begin{document}
\title{Non-commutative ergodic Theorems for actions of the hyperbolic groups}
\author{Genady Ya. Grabarnik\thanks{St. Johns University},
Alexander A. Katz\thanks{St. Johns University},
Laura Shwartz\thanks{IBM TJ Watson Research}}
\date{}
\maketitle

\begin{abstract}
The goal of this notice is to establish Not-commutative Pointwise Ergodic Theorems for actions of the Hyperbolic Groups. Similar non-commutative results were done by Bufetov, Khristoforov and Klimenko, and later by Pollicott and Sharp. We were interested to expand short notice in Policott and Sharp's paper about non-commutative ergodic theorems.

\end{abstract}

\section{Introduction}
After Birkhoff pointwise ergodic theorem, that may be considered as ergodic theorem for action of ${\mathbf Z}$, a number of authors worked on ergodic theorems for general groups. Among others we would like to mention work by Oceledetz \cite{Oseled}, Tempelman \cite{Tempelman}, Arnold and Krylov\cite{ArnKr}, Guivarch \cite{Guivarch}, Grigorchuk\cite{Grig86, Grig99}, Nevo and Stein\cite{Nevo94, NeSt}, Bufetov\cite{Buf99, Buf02, BufKK}, see also \cite{calegari2010combable} . Non commutative analogs are due to work of authors \cite{GKS, GKSI}, Anantharaman-Delaroche\cite{Anan}. Authors thank Danny Calegari for clarification of proper references and statement of Proposition \ref{thm:prop1}.   
\subsection{Settings}
Let $M$ be a von Neumann algebra with normal faithful state $\phi$, and $Aut(M)$ be a group of automorphisms of $M$ leaving state       $\phi$ invariant. 
We suppose that group $\Gamma$ acts on the $(M, \phi)$ as automorphisms preserving state $\phi$ and hence as  positive double stochastic operators, 
\begin{equation}
\gamma\cdot\phi=\phi;\  
\gamma\mathbbm{1}=\mathbbm{1} \text{, here } \gamma\in \Gamma
\end{equation}
We suppose that the action group $\Gamma$ is finitely generated group with symmetric set of generators $\Gamma_0$:
\begin{equation}
\Gamma={\cal{G}}(\Gamma_0);\ |\Gamma_0|<\infty;\ \gamma \in \Gamma_0, \text{ hence } \gamma^{-1}\in\Gamma_0.
\end{equation}
There is a natural notion of {\it length} on the group $\Gamma$ defined by choice of the generator set $\Gamma_0$
\begin{equation}
|\gamma |=inf\{n\geq 1| \gamma=\gamma_n\cdots\gamma_1,\ \gamma_k\in\Gamma_0-\{e\}\},
\end{equation}
here $\gamma\in \Gamma$  and  $e$  is a unit of $\Gamma$.

Based on the notion of length $|\cdot |$ we can consider metric $d$ on $\Gamma$ introduced as
\begin{equation}
d(\gamma, \delta)=|\gamma^{-1}\delta |,
\end{equation}
here $\gamma,\delta\in \Gamma$.

We suppose that the group we consider in the note satisfies introduced by Gromov hyperbolic property or the group is {\it word hyperbolic group},  for reference see \cite{Gromov}, meaning that there exists $\epsilon\geq 0$ such that for every geodesic triangle every point on the side of the triangle is within $\epsilon$ of both other sides, or triangle is $\epsilon$-thin.

We also suppose that group $\Gamma$ does not have cyclic subgroups of finite index. We will call such group $G$ not finitely cyclic (nfc) or non-elementary  (following Cannon \cite{Cannon}).

Let $N(n)$ be a cardinality of a set of all words of length n, and $\rho(G)$ (or just $\rho$) be an exponential growth rate of the group $\Gamma$:
\begin{equation}
N(n)=\#\{\gamma\in\Gamma | |\gamma |=n\},\ \rho(\Gamma)=\lim_{n\to\infty}N(n)^{1/n}.
\end{equation}

For every word hyperbolic group $\Gamma$ the limit $\rho(\Gamma)$ does exists (see for example Coornaert \cite{Coornaert}), and $N(n)$ is $\Theta(\rho^n)$, following Landau notations (for some positive constants $k_1, k_2$ holds $k_1\cdot \rho^n\leq N(n)\leq k_2\cdot \rho^n$).

Remind that sequence $X_n$ of operators from $M$ converges (bilaterally) almost sure to $X_0\in M$ if there exists sequence of orthogonal projections $E_k$ from $M$ increasing to $\mathbbm{1}$ such that 
\begin{equation}
\|(E_k)\cdot (X_n-X_0)\|_\infty\to 0\ (\text{or } \|E_k\cdot (X_n-X_0)\cdot E_k\|_\infty\to 0 ), \text{ when } k\to \infty,
\end{equation}
here $\|\cdot\|_\infty$ is a norm in von Neumann algebra $M$.

Now we are able to formulate the main result of the note (compare with \cite{PolShar}, Theorem 1)
\begin{theorem}
\label{thm:1}
Let $\Gamma$ be a word hyperbolic, nfc group with symmetric set of generators $\Gamma_0$. Let $\rho$ be an exponential growth constant for group $\Gamma$. Suppose that group $\Gamma$ acts on von Neumann algebra $M$ with faithful normal state $\phi$ as automorphisms, having state $\phi$ invariant.
Then for any $X\in M$, Cesaro averages of spherical averages
$$
s_n(X) =\frac{1}{N}\sum_{n=1}^N \rho^{-n}\sum_{|\gamma|=n+1} \gamma(X)
$$
converges in norm $L_1$ to some $\bar{X}$. 

Moreover $s_n(X)$ converges to $\bar{X}$ (bilaterally) almost sure.
\end{theorem}
\begin{remark}
Some generalization of the theorem \ref{thm:1} in the direction of more general groups (semigroups),
unboundedness of the operator $X$ with $X\in L_2(M, \phi)$ or $L_p(M, \tau)$ for$\tau$ being trace, etc.
may be found in the forthcoming work \cite{GA}.
\end{remark}

\section{Covering  Markov Group}
Let $\Gamma$ be a hyperbolic group and let $\Gamma_0$ be a symmetric generating set.
The action of the $\Gamma$ is associated till certain degree with subshift of finite type.

More specifically, consider finite directed graph $G$ with vertices $\Vertex G$ or just $\Ver$ and edges $\Edge G$ or just $\Edg$.
For the vertices $v_1,v_2\in \Ver$ let $e(v_1, v_2)$ be connecting them edge $e(v_1, v_2)\in\Edg$. Denote $v_1=I(e(v_1,v_2))$ and $v_2=F(e(v_1,v_2))$. We suppose that  $G$ has a vertex $v_0$ such that no edge ends in $v_0$. 

Let $\Pat(G)$ be the set of finite paths in digraph $G$ starting at $v_0$ ,
\begin{equation}
\Pat(G)=\{l=e_1e_2\cdots e_k \text{ with } I(e_j)=F(e_{j-1}), j=2,k\}.
\end{equation}
Denote by $|l|$ the length of the path $l$.

With digraph $G$ we associate transition matrix $A$ indexed by vertices $v\in\Ver$. 
We set element $a_{ij}$ of matrix $A$ to be equal to count of the edges with $I(e)=v_i, F(e)=v_j$:
$$
a_{ij}=\#\{e\in\Edg | I(e)=v_i, F(e)=v_j
$$

Remind that a square matrix $A$ with non-negative elements is irreducible if for any indices $i,j$ there is a power $n$ such that $A^n$ has $a_{ij}^{(n)}>0$ and reducible otherwise. If $n$ in the definition of irreducibility can be found uniformly by all $i,j$, matrix is aperiodic. 

Each edge $e\in \Edg(G)$ we associate (label) with automorphism $T_e$ of the algebra $M$. With every  path
$l=e_1e_2\cdots e_k$ we associate automorphism 
$$T_l=T_{e_1e_2\cdots e_k}=T_{e_1}\cdots T_{e_k}$$
that leaves state $\phi$ invariant. 
\begin{definition}
The hyperbolic word group $\Gamma$ and (any)  symmetric finite set of generators $\Gamma_0$ acts on algebra $M$ in a (strongly) Markov manner (or just (strong) Markov) if mapping of $\Gamma_0 \sim G_0 \mapsto Aut(M)$ may be lifted to injection 
$$\Gamma \sim \Pat(G)\mapsto Aut(M),$$
satisfying $|g|=|l|$ where $G\ni g=g_1\cdots g_k\sim l_{v_1\cdots v_k}\mapsto T_{v_1}\cdots T_{v_k}\in Aut(M). $
\end{definition}
\begin{note}
Since set of generators $\Gamma_0\sim G_0$ is symmetric, to make sense of injection statement we consider only paths that do not have two consecutive elements to be inverse one to another.
\end{note}
Since we associated group $\Gamma$ with the graph $G$, it will be not more than one edge connecting vertices $v_1,v_2\in\Ver$, and elements of the transition matrix $A$ are either 0 or 1.

\begin{figure}[ht]
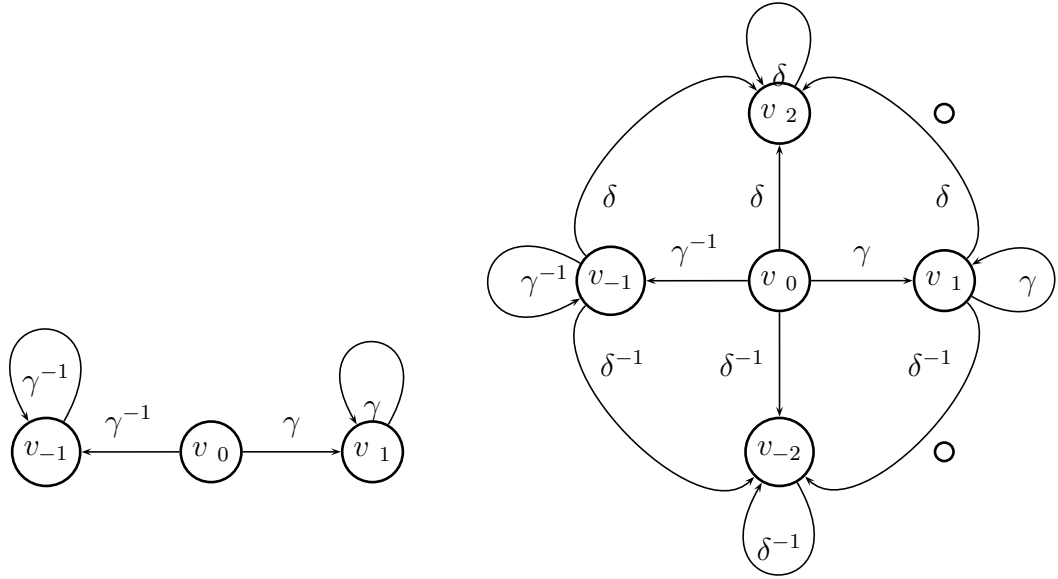

\begin{center}
		\psset{arrows=->,mnode=circle,linewidth=1pt}
		\begin{psmatrix}[colsep=1.3cm, rowsep=1.3cm]
		[name=s0]$v_{-1}$ & [name=s1]$v_{\ 0}$ & [name=s2]$v_{\ 1}$ 
		\psset{arrows=->,mnode=circle,linewidth=0.6pt}
		\ncline{s1}{s0}^{$\gamma^{-1}$}
		\ncline{s1}{s2}^{$\gamma$}
		\nccurve[angleA=60,angleB=120,ncurv=8]{s0}{s0}^[npos=.4]{$\gamma^{-1}$}
		\nccurve[angleA=60,angleB=120,ncurv=8]{s2}{s2}^[npos=.4]{$\gamma$}
		\end{psmatrix}
	\hspace{20mm}
		\psset{arrows=->,mnode=circle,linewidth=1pt}
		\begin{psmatrix}[colsep=1.3cm, rowsep=1.3cm]
		&  [name=r1]$v_{\ 2}$ & \\ 
		[name=s0]$v_{-1}$ & [name=s1]$v_{\ 0}$ & [name=s2]$v_{\ 1}$ \\ 
		&  [name=r2]$v_{-2}$ &
		\psset{arrows=->,mnode=circle,linewidth=0.6pt}
		\ncline{s1}{s0}^{$\gamma^{-1}$}
		\ncline{s1}{r1}<{$\delta$}
		\ncline{s1}{s2}^{$\gamma$}
		\ncline{s1}{r2}<{$\delta^{-1}$}
		\ncarc[arcangle=-90]{s2}{r1}<{$\delta$}
		\ncarc[arcangle=90]{s0}{r1}>{$\delta$}
		\ncarc[arcangle=90]{s2}{r2}<{$\delta^{-1}$}
		\ncarc[arcangle=-90]{s0}{r2}>{$\delta^{-1}$}
		\nccurve[angleA=150,angleB=210,ncurv=8]{s0}{s0}<[npos=.35]{$\gamma^{-1}$}
		\nccurve[angleA=60,angleB=120,ncurv=8]{r1}{r1}^[npos=.4]{$\delta$}
		\nccurve[angleA=-30,angleB=30,ncurv=8]{s2}{s2}<[npos=.35]{$\gamma$}
		\nccurve[angleA=-60,angleB=-120,ncurv=8]{r2}{r2}\Bput{$\delta^{-1}$}
		\end{psmatrix}
	\vspace{10mm}
\end{center}
\caption{Markov graphs for a) $\mathbb{Z}$ and b) $\mathbb{Z}^2$.}
\label{fig:MG}
\end{figure}

Now we consider a few examples of the (strongly) Markov groups.
\begin{example}
\begin{description}
\item[$\mathbb{Z}$] Since $\mathbb{Z}$ is a commutative group, let $z_1, z_2, e$ be a symmetric set of generators. Then associated digraph will have only edges from origin to $z_1, z_2$ and from $z_1, z_2$ to itself. The Markov digraph for the $\mathbb{Z}$ is shown in figure \ref{fig:MG} a).
\item[$\mathbb{Z}^n$] for $n\geq2 $ since the group is commutative, we chose linear order on generators and make sure that according to Markov digraph first we apply generators according to chosen order. The Markov digraph for the $\mathbb{Z}^2$ is shown in figure \ref{fig:MG} b).
\item[$\mathbb{F}_n$] To build Markov graph for the free group of $n$ generators we connect each vertex with every vertex except for it's inverse and origin.
\item[Word hyperbolic group] Any word hyperbolic group is strongly Markov (see for example \cite{GhysDelaharpe}).
\end{description}
\end{example}

 To translate nfc property into the transition matrix for digraph language it means that matrix at least does not have periodic irreducible component (with period higher than 1). Note that the first two examples, $\mathbb{Z}$ and $\mathbb{Z}^2$ have cyclic subgroups of finite index.

We now combine multiple facts either from textbooks or different sources into one proposition.
\begin{proposition}(cmp. \cite{PolShar}, \cite{calegari2010combable})
\label{thm:prop1}
Let $\Gamma$ be a word hyperbolic nfc group with symmetric generators set $\Gamma_0$. 
Suppose that $G$ is Markov graph of the group and $A$ is its transition matrix . 
By \cite{Coornaert} there exists $\rho$ with $C_1\rho^n\leq N(n)\leq C_2\rho^n$. Then 
\begin{enumerate}
\renewcommand{\theenumi}{\alph{enumi}}
\renewcommand{\labelenumi}{(\theenumi)}
\item\label{en:prop:1} 
The matrix $A$ is block lower triangular matrix having for
	$$
	\left( \begin{array}{cccc}
A_1 & 0 & 0 & 0 \\
B_{21} & A_2 & 0 &0 \\
B_{31} & B_{32} & A_3 & 0 \\
B_{l1} & B_{l2} & * & A_l 
\end{array} \right)
$$, where each $A_j$ is irreducible matrix.
\item\label{en:prop:2}  Spectral radius of the $A_{i_j}$ is $\rho$, for some $i_j$, and for the others $i$ spectral radius of $\rho(A_j)\leq\rho_1<\rho$. We call such $A_{i_j}$ contributing.
\item\label{en:prop:3}  For any $i_1,i_2\in 1,l$ being indices of the contributing components
\begin{equation}
B_{i_1, i_{j_1}}\cdots B_{i_{j_n}, i_2}=0; \qquad A_{i_1}^mB_{i_1, i_2}A_{i_2}^n=0,
\end{equation}
in other words any path intersect contributing components only once (but may loop in it).
\item\label{en:prop:4}  
Number of paths outside of contributing components grow asymptotically not faster than  $\rho_1^n$.
\end{enumerate}
\begin{proof}

The statement \ref{en:prop:1}) is a standard statement about Markov matrices, with additional fact that for word hyperbolic nfc group states of a transition matrix are connected 

The statement \ref{en:prop:2}) is also a standard statement about Markov matrices. 

The statement \ref{en:prop:3}) is due to \cite{calegari2010combable}; the first part follows from the Coornaert's result stated in the formulation of the proposition and 
consideration of the $A^m$ having addends $A_{i_1}^lB_{i_1, i_{j_1}}\cdots B_{i_{j_n}, i_2}A_{i_2}^{m-n-2}$ for every $1<n<m-2$. 
Since $A_{i_j}^l$ elements grow as $\rho^l$ and are irreducible, 
that means that if $B_{i_1, i_{j_1}}\cdots B_{i_{j_n}, i_2}$ is not zero we will have a number  proportional to $n$ of elements 
proportional to $\rho^n$, which contradict Coornaert's estimate. The statement \ref{en:prop:3}), the second part is treated in similar manner.

The statement \ref{en:prop:4}) follows from consideration of $A^n$, \ref{en:prop:3} of the Proposition, the fact that there may be product of powers of $A_i$ and $B_{ij}$ not longer than $2l-1$, and part \ref{en:prop:2}) of the Proposition.
\end{proof}
\end{proposition}

To work with the group action we consider construction of the covering Markov group 
which is due to Grigorchuk \cite{Grig99}, J.-P. Thouvenot (oral communication), and Bufetov \cite{Buf99}.

Consider $M_k=M\times \mathbb{Z}_k$, 
and define covering Markov operator $P_{\Gamma, \Gamma_0, G}$ depending on the symmetric finite subset $\Gamma_0$ (with $\#(\Gamma_0)=k$) 
of the generators of group $\Gamma$ and Markov graph $G$ of the group $\Gamma$ as
\begin{equation}
(P_{\Gamma, \Gamma_0, G} (x_i, i)_{i=1,k})_m=\sum_{j=1,k}A_{m,j}( \gamma_{v_m,v_j}x_m,j), 
\label{eq:coveringMarkov}
\end{equation}
here $(x_i,i)_{i=1,k}\in M_k$. We will omit parameters of operator $P_{\Gamma, \Gamma_0, G}$ and use just $P$.
From the definition of the covering Markov operator it follows that for $M_k\ni \bar{x}=(x, i)_{i=1,k},\ x\in M$
\begin{equation}
\sum_{\gamma\in\Gamma, |\gamma |=n}\gamma(x)=P^n(x).
\label{eq:covMarkPath}
\end{equation}
To use special structure of the nfc word hyperbolic groups we partition operator P on 
action of block diagonal part of matrix $A$ (operator $D$) and sub block diagonal part (operator $Q$):
\begin{equation}
(D (x_i, i)_{i=1,k})_m= \sum_{j, m\in\text{ block }p, all p}A_p(m,j)( \gamma_{v_m,v_j}x_m,j), 
\label{eq:D}
\end{equation}
\begin{equation}
(Q (x_i, i)_{i=1,k})_m= \sum_{j, m \text{ in different blocks }p_1,p_2}B_{p_1,p_2}{m,j}( \gamma_{v_m,v_j}x_m,j), 
\label{eq:Q}
\end{equation}
The expression \ref{eq:covMarkPath} may be simplified by taking into account parts \ref{en:prop:3}) and  \ref{en:prop:4})  of Proposition \ref{thm:prop1}:
\begin{equation}
\sum_{\gamma\in\Gamma, |\gamma |=n}\gamma(x)=\sum_{k=1,m} 
\sum_{\substack{ 1\leq i_1<\cdots<i_k\leq k, \\  l_i>=0, \sum_{j=1}^kl_j=n-k+1}} A_{i_1}^{l_1}B_{i_1i_2}A_{i_2}^{l_2}B_{i_2i_3}\cdots B_{i_{k-1}i_k}A_{i_{k}}^{l_k},
\label{eq:markPath}
\end{equation} 
and in each product can have not more than one $A_j$ is contributing.

We extend state $\phi$ on $M_k$ as $\phi_k((x_j)_{j=1,k})=1/k\sum_{j=1,k}\phi(x_j)$.

In order to use the construction above, we need to normalize operator $P$ in the following manner.
Due to Peron-Frobenius theorem, there is an eigenvector $w_j$ with 
all positive coordinates for matrices $A_j$ such that $A_jw_j=\rho w_j$.
We define stochastic matrix $R(i,j)=\rho^{-1}A_l(i,j) w(j)/w(i)$ for $i,j\in supp(A_j)$, 
and its action on $M_k, \phi_k$ for $(x_i)\in M_k$ as
$$
\mathcal{P}(x)(i)=\sum_{j=1,k} R(i,j)\gamma_{i,j}((x_j), j).
$$
The operator $\mathcal{P}$ is Markov operator. 
Let $L_w$ be the operator on $M_k$ acting by multiplication on $w_j$.
Due to definition of the matrix $R(i,j)$ from the chain rule, it follows that
\begin{equation}
D^q=\rho^q L_w R (L_w)^{-1}.
\label{eq:D1}
\end{equation}

Consider now 
\begin{equation}
A_N(x)=\frac{1}{N}\sum_{j=0}^N\rho^{-j} \sum_{|\gamma|=j+1} \gamma(x).
\label{eq:AN}
\end{equation}
By previous consideration,
\begin{align}
A_N(x)=\frac{1}{N}\sum_{n=2}^N \rho^{-n}\sum_{k=1,m} 
\sum_{\substack{ 1\leq i_1<\cdots<i_k\leq k, \\  l_i>=0, \sum_{j=1}^kl_j=n-k+1}} & \rho^{l_1} L_w R^{l_1} (L_w)^{-1}(\rho^{-1}Q)\times \\
& (\rho^{-1}Q)\rho^{l_k} L_w R^{l_k} (L_w)^{-1}
\label{eq:AND}
\end{align}

We now are ready to establish pointwise ergodic theorems in von Neumann algebras.
\section{Pointwise Ergodic Theorems in von Neumann Algebras}
First we reduce arbitrary $x\in M$ to positive  $x\in M_+$ and establish non commutative majorant ergodic theorem for action of the  nfc word hyperbolic groups (see for example \cite{FuNe} for commutative analog). The proof follows 'standard' schema for proof of ergodic theorems. We establish majorant type theorem, then show  mean ergodic theorem (in $|| .||_2$). 

\begin{theorem}[Majorant Ergodic theorem]
\label{thm:maj}
Let $x\in M_+$ be a positive operator from $M$.
 Then there exists $\tilde{x}\in M_+$ with $||\tilde{x}||\leq C(G,||x||) $ with $C(G,||x||)\to 0$ for $||x||\to 0$, such that 
\begin{equation}
A_n(x) \leq \tilde{x}.
\end{equation}
\end{theorem}
\begin{proof}
The proof of the theorem is consequence of the \ref{en:prop:4})  of Proposition \ref{thm:prop1}.
We can sum up to $y\in M_+$ all addends of the equation \ref{eq:AND} without contributing component,
 it will converge  since its norm is majorated by constant times decreasing geometric progression ( $C_2(G)||x||$).

Next we partition the sum by different contributing components and to each contributing component  apply Lance's  majorant theorem  for the  Cesaro  to get majorant $z\in M_+$ with norm of scale of square root of norm of $y$.

Next we apply again sum of not contributing components and again estimate norm of majorant as $\tilde{x}$ as $C_2(G)||z||$.
\end{proof}

\begin{theorem}[Mean Ergodic theorem]
Let $||x||_2$ be a $\phi( x^*x)^\frac{1}{2}$ for $x\in M$. Then for some $x_0\in M$
\begin{equation}
A_n(x)\xrightarrow{||.||_2}x_0
\end{equation}
\end{theorem} 
\begin{proof}
For the fixed $\epsilon$ we find such a big $m$ that norm of the all non contributing components with total power not less than $m$ does not exceed $\epsilon/3$. This follows from the \ref{en:prop:4})  of Proposition \ref{thm:prop1}. 

By applying  mean ergodic theorem to each of the contributing components' Cesaro sums we will make norm of the difference  $A_n(x)-x_0 $ small.

By partitioning again sum of all not contributing components on total power less than $m$ we get finite number of addends and 
will control them by smallness resulted from mean ergodic theorem. For the total power of non contributing components greater than $m$ we use estimate $\epsilon/3$ (see above).

\end{proof}

To deal with (double side) almost everywhere convergence 
we need to consider (o)-convergence. Remind that $(o)-$ convergence considered in the \cite{GoldGra}.

The sequence of self adjoin operators $\{x_n\}_{n=1}^\infty\subset M_{sa}$ 
(o) converges to the $x_0\in M_{sa}$ 
if there exists monotonically decreasing to $0$ sequence $\{y_n\}_{n=1}^\infty\subset M_+$
such that 
\begin{equation}
-y_n\leq x_n-x_0\leq y_n.
\label{eq:0converges}
\end{equation}

\begin{theorem}[(o) convergence]
\begin{equation}
A_n(x)\xrightarrow{(o)}x_0
\label{eq:0convergence}
\end{equation}
\end{theorem}
\begin{proof}
Since $A_n(x)\xrightarrow{||.||_2}x_0$ we chose subsequence in $n_j$ with both difference  $||A_n(x)-x_0||_2$ and $C(G,||A_n(x)||_2)$ decreasing  than $2^{-j}$. 

Again using   the \ref{en:prop:4})  of Proposition \ref{thm:prop1} we conclude that 
\begin{equation}
||A_n(x)-A_n(A_{n_j}(x)||<\frac{C_3}{n}||x||
\label{eq:0convergence}
\end{equation}
in norm of $M$ for $n>n_j'=max{n_j, 2^{2j}}$
By summing up all majorants for all $n_j'$ starting from $j_0$ and adding scaled constants $2^{-2j}\ONE$ 
we get decreasing sequence of $y_{0j}$ of positive operators in $L_2(M, \phi)$.
\end{proof}
Due to Schwartz type inequality we can majorate (for positive $x$) $(A_N(x))^2$ 
by scale of the $A_N(x^2)$, this consideration also implies (o) convergence of $(A_N(x))^2$. 

Since t (o) convergence imply double side almost everywhere convergence, 
and (o) convergence of squares imply almost everywhere convergence.
This establishes theorem \ref{thm:1}

\section{Extensions}
In this section we outline different possible extensions of the Theorem \ref{thm:1}.

Since we did not used the fact that $\Gamma$ is group an only used Proposition \ref{thm:prop1},
so theorem \ref{thm:1} is valid for the semigroup of positive double contractions having structure 
satisfying results of Proposition \ref{thm:prop1}.

Next we can extend results of the Theorem \ref{thm:1} to the non-commutative 
(Haagerup's) $L_p$ spaces ($p>1$) and action of the commuting with modular operator positive double contraction.
We use the existence of  majorant based on Theorem 7.9 of the Haagerup, Junge, Xu work \cite{HJX} on reduction method. 
This will give $||.||_p$ and double side almost everywhere convergence.
Again, for $p>2$, we use Schwartz type inequality to obtain majorant for squares of $A_N$, 
and obtain theorem for the almost everywhere convergence.

Extensions to the Jordan algebras is considered in the forthcoming \cite{GA}.

\end{document}